\documentclass[12pt,leqno]{article}
\usepackage{amssymb,amsmath,amsfonts,amsbsy, xspace, latexsym}

\newtheorem{theorem}{Theorem}[section]
\newtheorem{lemma}{Lemma}[section]
\newtheorem{proposition}{Proposition}[section]
\newtheorem{definition}{Definition}[section]
\newtheorem{corollary}{Corollary}[section]
\newtheorem{notation}{Notation}[section]
\newtheorem{simpler notation}{Simpler Notation}[section]
\newtheorem{remark}{Remark}[section]
\newtheorem{example}{Example}[section]

\newcommand{\Proof}{\noindent \textbf{Proof:\;}}

\newcommand{\bra}{\ensuremath{\large\langle}}

\newcommand{\ket}{\ensuremath{\large\rangle}}

\newcommand{\card}{\ensuremath{{\rm{card}}}}

\newcommand{\st}{\ensuremath{{{\rm{st}}}}}

\newcommand{\cl}{\ensuremath{{\rm{cl}}}}

\begin{document}

\title{HAHN FIELD REPRESENTATION OF\\
 A. ROBINSON'S ASYMPTOTIC NUMBERS}

\author{Todor D. Todorov (ttodotrov@calpoly.edu)\\
			Mathematics Department\\		
			California Polytechnic State University\\
		 San Luis Obispo, CA 93407, USA
						\and
			Robert S. Wolf (rswolf@calpoly.edu)\\
			Mathematics Department\\		
			California Polytechnic State University\\
			San Luis Obispo, California 93407, USA}
\date{}
\maketitle	

\begin{abstract}
Let $^*\mathbb{R}$ be a nonstandard extension of $\mathbb{R}$ and $\rho$ be a positive infinitesimal in
$^*\mathbb{R}$. We show how to create a variety of isomorphisms between A. Robinson's
field of asymptotic numbers $^\rho\mathbb{R}$ and the Hahn field 
$\widehat{^\rho\mathbb{R}}(t^\mathbb{R})$, where $\widehat{^\rho\mathbb{R}}$
is the residue class field of $^\rho\mathbb{R}$. Then, assuming that $^*\mathbb{R}$ is fully saturated we
show that
$\widehat{^\rho\mathbb{R}}$ is isomorphic to $^*\mathbb{R}$ and so $^\rho\mathbb{R}$ contains a copy of
$^*\mathbb{R}$. As a consequence (that is important   for applications in non-linear theory of
generalized functions) we show that every two fields of asymptotic numbers corresponding to different
scales are isomorphic.
\end{abstract}
\section{Introduction} \label{S: Introduction} The main purpose of this
article is to study in detail some properties of Abraham Robinson's
field of {\em asymptotic numbers} $^\rho\mathbb{R}$ (known also as {\em Robinson's
valuation field}) and its residue class field $\widehat{^\rho\mathbb{R}}$ (A. Robinson
(\cite{aRob73} and A.H. Lightstone and A. Robinson~\cite{LiRob}). In this sense the
article is a continuation of the handful of papers written on the subject
(known as {\em non-standard asymptotic analysis}), including:
W.A.J. Luxemburg~\cite{wLux}, Li Bang-He~\cite{liBang}, B. Diarra~\cite{bDi} and
V. Pestov~\cite{vPes}. Our emphasis, however, is on those properties of
$^\rho\mathbb{R}$ and $\widehat{^\rho\mathbb{R}}$ which seem to be of importance for some more recent
applications of $^\rho\mathbb{R}$ to the non-linear theory of generalized functions
and partial differential equations: M. Oberguggenberger~\cite{mOber88}-\cite{mOber95}, M. Oberguggenberger
and T. Todorov~\cite{OberTod}, T. Todorov~\cite{tTod88}-\cite{tTod99}. In particular, 
the complex version
$^\rho\mathbb{C}$ of the field $^\rho\mathbb{R}$ could be viewed as a ``non-standard counterpart'' of
Colombeau's ring of generalized numbers $\overline{\mathbb{C}}$ (J.F. Colombeau~\cite{jfCol85}).
Similarly, the algebra of asymptotic functions $^\rho\mathcal{E}(\Omega)$ (M. Oberguggenberger and T.
Todorov~\cite{OberTod})  can be viewed as a ``non-standard counterpart'' of a typical algebra of
generalized functions
$\mathcal{G}(\Omega)$ in Colombeau's theory, where the constant functions in $^\rho\mathcal{E}(\Omega)$
are  exactly the
asymptotic numbers in $^\rho\mathbb{C}$ (T. Todorov~\cite{tTod99}). We should mention that
$^\rho\mathbb{C}$ is an algebraically closed Cantor complete non-Archimedean field while
$\overline{\mathbb{C}}$ is a non-Archimedean ring with zero divisors. Thus, the involvement of the
non-standard analysis results to an essential improvement of the properties of the generalized scalars,
compared with the original Colombeau's theory of generalized functions.

\section{Preliminaries: Non-Archimedean Fields} \label{S: prelim}

In this section we cover various preliminary material. Given a
(totally) ordered ring, we define \emph{absolute value} in the usual way: $|x| =
\max\{x, -x\}$. \emph{Intervals} and the \emph{order topology} are also
defined as usual. Unless specified otherwise, we let every ordered
structure be topologized with the order topology. We assume the reader
is familiar with the basics of \emph{non-archimedean fields}. If $\mathbb{K}$ is a
nonarchimedean field and $x, y \in \mathbb{K}$, we write $x \approx y$ to mean
that $x-y$ is infinitesimal. For more details, see T. Lightstone and
A. Robinson~\cite{LiRob}.
	
In any ordered field, we define the usual notions of \emph{fundamental} or
\emph{Cauchy} sequences and \emph{convergent} sequences, with respect to
absolute value. In a non-archimedean field, a sequence $(a_n)$ is fundamental if
and only if the sequence of differences $(a_{n+1}-a_n)$ converges to
$0$. This fact simplifies the theory of series in a sequentially (Cauchy)
complete non-archimedean field: a series converges iff the terms
approach 0, and there is no such thing as conditional convergence.

\begin{definition} \label{D:Cantor} {\em (Cantor Completeness)}: Let $\kappa$ be an infinite
cardinal. An ordered field is called \textbf{$\kappa$-Cantor complete} if every collection of
fewer than $\kappa$ closed intervals with the finite intersection property
(F.I.P.) has nonempty intersection. An ordered field is simply called {\bf
Cantor complete} if it is $\aleph_1$-Cantor complete (where $\aleph_1$ denotes the successor of
$\aleph_0=\card(\mathbb{N})$). This means every nested sequence of closed intervals has
nonempty intersection.
\end{definition}

	It is easy to show that a Cantor complete ordered field must be
sequentially complete. Two counterexamples to the converse are described
in the discussion of Hahn fields later in this section. 

Recall that the usual definition of $\kappa$-saturation in non-standard
analysis is obtained by replacing ``closed interval'' by ``internal
set'' in the above definition. There is also a model-theoretic or
first-order definition of $\kappa$-saturation, in which ``internal set''
is replaced by ``definable set'' (and which for {\em real-closed} fields is
equivalent to replacing ``closed interval'' by ``open interval''). This
is weaker than the non-standard notion, but still much stronger than
$\kappa$-Cantor completeness. 

\begin{example}[Non-Standard Real Numbers]\label{Ex: Nonstandard Real Numbers}  A typical examples for
non-archimedean fields are the non-standard extensions of $\mathbb{R}$ in the framework of the
non-standard analysis : A set $^*\mathbb{R}$ is called a \textbf {non-standard extension} of $\mathbb{R}$ 
if $\mathbb{R}\subsetneqq{^*\mathbb{R}}$ and there exists a mapping $*$ from
$V(\mathbb{R})$ into {\em $V(^*\mathbb{R})$} which maps $\mathbb{R}$ at $^*\mathbb{R}$, 
and which satisfies
the Transfer Principle (A. Robinson~\cite{aRob66}). Here $V(\mathbb{R})$ and {\em $V(^*\mathbb{R})$}
denote the superstructures on $\mathbb{R}$ and $^*\mathbb{R}$, respectively. Thus, by  definition, any
non-standard extension $^*\mathbb{R}$ of $\mathbb{R}$ is a totally ordered real closed field 
extension of $\mathbb{R}$ (with respect to the operations inherited from $\mathbb{R}$), since 
$\mathbb{R}$ is a totally ordered real closed field.  Also, $^*\mathbb{R}$ must be a 
\textbf{non-archimedean field} as a proper field extension of
$\mathbb{R}$. For more details on non-standard analysis we refer to A. Robinson~\cite{aRob66} and T.
Lindstr\o m~\cite{tLin}, where the reader will find many other references to the subject. For a
really brief introduction to the subject - see the Appendix in T. Todorov~\cite{tTod96}. 
\end{example}

\begin{definition}[Valuation] \label{D: Valuation} Let $\mathbb{K}$ be an ordered field. A function
$v$ from $\mathbb{K}$ to $\mathbb{R} \cup \{\infty\}$ is called a \textbf {valuation} on $\mathbb{K}$
provided (for any $x,y \in \mathbb{K}$): 

	(a) $v(x) = \infty$ if and only if $x = 0$ 

	(b) $v(xy) = v(x) + v(y)$ (\textbf{logarithmic property})

	(c) $v(x+y) \ge \min\{v(x),v(y)\}$ (\textbf{non-archimedean property})

	(d) $|x|<|y| $ implies $v(x) \ge v(y)$ (\textbf{convexity} or \textbf{
compatibility with the ordering on} $\mathbb{K}$). If $v$ is a valuation on $\mathbb{K}$, the pair 
$(\mathbb{K}, v)$ is called a \textbf{valuation field}. (Here it is understood that $x \le \infty$ and 
$x+\infty = \infty$, for every $x$ in $\mathbb{R} \cup \{\infty\}$.)
\end{definition}
	
	We give $\mathbb{R}\cup\{\infty\}$ the usual order topology: all ordinary open
intervals and all intervals of the form $(a, \infty]$ are basic open
sets. If $v$ is a valuation then it is easy to show that $v(1) = 0, v(-x) =
v(x), v(1/x) = -v(x)$ whenever $x \ne 0$, and $v(x) \ne v(y)$ implies
$v(x+y) = \min\{v(x),v(y)\}$. For additional details on valuations,
we refer the reader to P. Ribenboim \cite{pRib} or A.H. Lightstone and A. Robinson
\cite{LiRob}.

The {\bf trivial valuation} on $\mathbb{K}$ is the one defined by $v(x) = 0$ for
every non-zero $x\in\mathbb{K}$. This is the only possible valuation on an
archimedean field.

Given a valuation field $(\mathbb{K}, v)$, we can define several important
structures: $\mathcal{R}_v(\mathbb{K}) = \{ x \in \mathbb{K}\mid\,v(x) \ge 0\,\}$ is a convex 
subring of
$\mathbb{K}$ called the {\bf valuation ring} of $(\mathbb{K}, v)$. $\mathcal{I}_v(\mathbb{K}) = 
\{\,x \in
\mathbb{K} \mid v(x) > 0 \,\}$ is also convex and is the unique maximal ideal in
$\mathcal{R}_v(\mathbb{K})$, called the {\bf valuation ideal} of $(\mathbb{K}, v)$. And
$\widehat{\mathbb{K}}$ =
$\mathcal{R}_v(\mathbb{K})/\mathcal{I}_v(\mathbb{K})$ is an ordered field called the {\bf residue class
field} of $(\mathbb{K}, v)$. Also, we let $\mathcal{U}_v(\mathbb{K}) =\mathcal{R}_v(\mathbb{K})
\setminus\mathcal{I}_v(\mathbb{K}) = \{\,x
\in
\mathbb{K}\mid v(x)=0\,\}$, the multiplicative group of units of $\mathcal{R}_v(\mathbb{K})$. The elements
of the multiplicative group $\mathcal{I}_v(\mathbb{K}_+) = \{\,x \in
\mathbb{K}_+ \mid v(x) > 0 \,\}$ are the \textbf{scales} of $\mathbb{K}$.  Finally, the {\bf valuation
group} of $(\mathbb{K}, v)$, denoted
$G_v$, is the image of
$\mathbb{K}\setminus\{ 0\}$ under $v$, a subgroup of
$\mathbb{R}$.

A valuation $v$\, on a field $\mathbb{K}$ induces a metric $d$, called a {\bf
valuation metric}, by the rule $d(a,b) = e^{-v(a-b)}$, with the
understanding that $e^{- \infty} = 0$ This metric satisfies the {\bf
ultrametric inequality}: $d(a,c) \le \max \{ d(a,b), d(b,c) \} $, for
all $a,b,c \in \mathbb{K}$. The ultrametric inequality has some strange
consequences. For one thing, all triangles are isosceles - a triangle
can have a shortest side but not a longest one. Also, any two closed
balls are either disjoint or one is a subset of the other. Furthermore,
every element of a ball is a center of the same ball. That is, $B(a,r) =
B(c,r)$ whenever $c \in B(a,r) $, where the notation $B$ denotes either
open or closed balls.

The topology defined on $\mathbb{K}$ from its valuation metric is called the {\bf
valuation topology}. We can also define the valuation topology directly from $v$;
this is the topology with basic open sets: $ \{\, x \in \mathbb{K} \mid v(x-a) > n\,\} $, where $a \in
\mathbb{K},\, n\in \mathbb{N}$.

The following facts are straightforward:

\begin{proposition}[Valuation and Topology] \label{P: ValTop} Let $(\mathbb{K}, v)$ 
be a nontrivial
valuation field. Then:

(a) The valuation topology and the order topology on $\mathbb{K}$ are the same.

(b) The functions\, $v$ \,and\, $d$\, are continuous.

(c) The notions of fundamental sequences, convergent sequences, and
sequential completeness with respect to the valuation metric coincide
with those notions with respect to absolute value. 

(d) For any infinite sequence $(a_n)$ of elements of $\mathbb{K}$, $a_n\rightarrow 0$ if and only if
$v(a_n) \rightarrow \infty$. 
\end{proposition}
\begin{definition}[Spherically Complete] \label{D:SpherCompl} 
A metric space is called {\bf spherically complete} if every nested sequence of closed balls has
nonempty intersection.
\end{definition}

A spherically complete metric space must be sequentially complete. Also,
if a {\em metric space, generated by a valuation}, is spherically complete, then
every collection of closed balls with the F.I.P. has nonempty intersection.

In contrast to Proposition~\ref{P: ValTop}(c) above, the following cannot be reversed:

\begin{theorem} \label{T: CantorSpher}
If $\mathbb{K}$ is Cantor complete, then $\mathbb{K}$ is spherically complete
with respect to any valuation on it. 
\end{theorem}

\noindent{\bf Proof:} Assume $\mathbb{K}$ is Cantor complete and $v$ is a valuation
on $\mathbb{K}$. Let $(B_n)$ be a strictly decreasing sequence of closed balls in
$\mathbb{K}$. Our goal will be to define, for each $n$, a closed interval $I_n$
such that $B_{n+1} \subseteq I_n \subseteq B_n$. This will make $(I_n)$
a decreasing sequence of closed intervals, so by Cantor completeness,
$\cap I_n \ne \emptyset$. But any element of $\cap I_n$ is clearly also
in $\cap B_n$, establishing that $\mathbb{K}$ is spherically complete.

It remains to define $I_n$. Say $B_n = B(a_n , r_n)$, the closed ball of
radius $r_n$ centered at $a_n$. We have $a_{n+1} \in B_n$, so by our
earlier remarks, $a_{n+1}$ is also a center of $B_n$. That is, $B_n =
B(a_{n+1} , r_n)$. Since $B_{n+1} \subset B_n$, we can choose some $c
\in B_n - B_{n+1}$. Thus, $r_{n+1} < d(a_{n+1}, c) \le r_n$. Let $s =
|a_{n+1} - c|$. If $c < a_{n+1}$, let $I_n = [c, a_{n+1} + s]$. If $c >
a_{n+1}$, let $I_n = [a_{n+1} - s, c]$. Since the distance between two
numbers depends only on their difference, we have $B_{n+1}\subseteq I_n \subseteq B_n$, 
as desired. $\blacktriangle$ 

	Thus, for a valuation field, $\aleph_1$-saturation implies Cantor
completeness implies spherical completeness implies sequential
completeness. A simple counterexample to the first converse is $\mathbb{R}$.
Counterexamples to the latter two converses are described below.

	An important category of non-archimedean fields are {\bf generalized
power series fields} or {\bf Hahn fields}, introduced in H. Hahn \cite{hH}.
Let $\mathbb{K}$ be a field and $G$ an ordered abelian group. (In this paper, $G$ is
usually a subgroup of $(\mathbb{R}, +))$. For any formal power series $f =
\sum_{g \in G} {{a_g} t^g}$, where each $a_g \in \mathbb{K}$, the \textbf{support}
of $f$ is $\{\,g \in G \mid a_g \neq 0\,\}$. Then the set of all such $f$'s
whose support is well-ordered is a field (using ordinary polynomial-like
addition and multiplication) denoted by $\mathbb{K}(t^G)$ or $\mathbb{K}((G))$. $\mathbb{K}$ is
naturally imbedded in $\mathbb{K}(t^G)$ by mapping any $a$ in $\mathbb{K}$ to $at^0$.

$\mathbb{K}(t^G)$ has a canonical $G$-valued Krull valuation in which each
nonzero power series is mapped to the least exponent in its support. If
$\mathbb{K}$ is ordered, then $\mathbb{K}(t^G)$ has a natural ordering in which an element
is positive if and only if the coefficient corresponding to the least
element in its support is positive. This ordering is compatible with the
canonical valuation, and is the unique ordering on $\mathbb{K}(t^G)$ in which
every positive power of $t$ is between 0 and every positive element of
$\mathbb{K}$.

One field of this type that has been studied extensively is
$\mathbb{R}(t^{\mathbb{Z}})$, which is called the field of \textbf {Laurent series}.
Another is the subfield $\mathbb{R}\bra t^{\mathbb{R}}\ket$ of $\mathbb{R}(t^{\mathbb{Z}})$ 
consisting of
those series whose support is either a finite set or an unbounded set of order type $\omega$. (We will
describe this field in a more concrete way in Section 3). This field was introduced by
Levi-Civita in \cite{tLC} and later was investigated by D. Laugwitz in
\cite{dLaug} as a potential framework for the rigorous foundation of
infinitesimal calculus before the advent of Robinson's non-standard
analysis. It is also an example of a real-closed valuation field that is
sequentially complete but not spherically complete (see V. Pestov 
\cite{vPes}, p. 67).
 
From W. Krull \cite{wKrull} and Theorem 2.12 of W.A.J. Luxemburg \cite{wLux}, it is known
that every Hahn field of the form $\mathbb{K}(t^\mathbb{R})$ is spherically complete in its canonical
valuation. In particular, $\mathbb{Q}(t^{\mathbb{R}})$ is spherically complete. But
$\mathbb{Q}(t^{\mathbb{R}})$ is not Cantor complete, for the same reason that $\mathbb{Q}$ is not Cantor
complete.

\section{Asymptotic and Logarithmic Fields} \label{S: Asy}

Throughout this paper, we let $^*\mathbb{R}$ be a non-standard extension of
$\mathbb{R}$ (Example~\ref{Ex: Nonstandard Real Numbers}). Every standard set, relation and function 
$X$ has a nonstandard extension $^*X$; we may omit the symbol $*$, when no ambiguity
arises from doing so. Until Section 5, we only need to assume that
$^*\mathbb{R}$ is $\aleph_1$-saturated.

	Let $\rho\in{^*\mathbb{R}}$ be a fixed positive infinitesimal. Following
A. Robinson \cite{aRob73} we define the field of {\bf Robinson's real
$\rho$-asymptotic numbers} as the factor space
$^\rho\mathbb{R}=\mathcal{M}_\rho(^*\mathbb{R})/\mathcal{N}_\rho(^*\mathbb{R})$,  where
\begin{align}\notag
&\mathcal{M}_\rho(^*\mathbb{R})=\{x\in{^*\mathbb{R}}\mid |x|\leq\rho^{-n}\text{\,for some\,} 
n\in\mathbb{N}\},\\ 
&\mathcal{N}_\rho(^*\mathbb{R})\,=\,\{x\in{^*\mathbb{R}}\mid |x|\leq\rho^{n}\text{\,for all\,} 
n\in\mathbb{N}\},\notag
\end{align}
\noindent are the sets of $\rho${\bf -moderate} and $\rho${\bf -null}
numbers in $^*\mathbb{R}$, respectively. It is easy to check that $\mathcal{M}_\rho(^*\mathbb{R})$ is
a convex subring of $^*\mathbb{R}$ and $\mathcal{N}_\rho(^*\mathbb{R})$ is a convex maximal ideal
of $\mathcal{M}_\rho(^*\mathbb{R})$. Therefore, $^\rho\mathbb{R}$ is a totally ordered field. The
quotient mapping $q:\mathcal{M}_\rho(^*\mathbb{R}) \to{^\rho\mathbb{R}}$, when restricted to
$\mathbb{R}$, becomes an embedding of $\mathbb{R}$ into $^\rho\mathbb{R}$. The asymptotic
number $q(\rho)$, denoted by $s$, or by $\overline{\rho}$, is called the \textbf{canonical scale} 
of $^\rho\mathbb{R}$.

It is simple to show that $^\rho\mathbb{R}$, like $^*\mathbb{R}$, is a real-closed,
Cantor complete, non-archimedean extension of $\mathbb{R}$. The intervals
$[-s^n, s^n]$, where $n \in \mathbb{N}$, form a base for the neighborhoods of
$0$ in the order topology on $^\rho\mathbb{R}$. So $^\rho\mathbb{R}$ is first countable, like
$\mathbb{R}$ but unlike $^*\mathbb{R}$. 

Let $f$ be any function whose domain and range are subsets of $\mathbb{R}$.
Chapter 4 of (A. Lightstone and A. Robinson \cite{LiRob}) describes the natural way
of using $^*f$ to try to define a function\, $^\rho\!f$\, that is an
extension of $f$, and whose domain and range are subsets of $^\rho\mathbb{R}$.
This procedure works for many but not all real functions. In particular,
it works for all polynomials. It also works for the logarithm function,
so we obtain a function\, $^\rho\!\ln: {^\rho\mathbb{R}_+} \to {^\rho\mathbb{R}}$
satisfying $^\rho\!\ln(xy) = {^\rho\!\ln(x) + ^\rho\!\ln(y)}$, etc. 
We shall write simply $\ln$ instead of $^\rho\!\ln$ when no confusion could arise. 
On the other hand, the exponential function cannot be defined in $^\rho\mathbb{R}$, since 
$e^{1/\rho} \notin \mathcal{M}_\rho(^*\mathbb{R})$.

There is a canonical valuation $v$ on $^\rho\mathbb{R}$ defined by $v(q(\alpha))
= \st(\log_{\rho}|\alpha|)$, where $\st$ is the standard part mapping in
$^*\mathbb{R}$ and $\alpha \in \mathcal{M}_\rho(^*\mathbb{R})$. It is easily shown that $v$ is
well-defined on $^\rho\mathbb{R}$ and is indeed a valuation (see A. Lightstone and
A. Robinson \cite{LiRob}, p. 79). In particular, we have $v(s)=1$. Let
$\mathcal{R}_v(^\rho\mathbb{R})=\{x\in {^\rho\mathbb{R}}\mid v(x)\geq 0  \}$ and
$\mathcal{I}_v(^\rho\mathbb{R})=\{x\in{^\rho\mathbb{R}}\mid v(x)> 0 \}$ be the valuation ring and the
valuation ideal of
$(^\rho\mathbb{R}, v)$, respectively, and
$\widehat{^\rho\mathbb{R}}=\mathcal{R}_v(^\rho\mathbb{R})/\mathcal{I}_v(^\rho\mathbb{R})$ 
be the residue class field of
$(^\rho\mathbb{R}, v)$ (called also {\em logarithmic field}). Notice that $\widehat{^\rho\mathbb{R}}$ can
be viewed as well  as a quotient field of a subring of
$^*\mathbb{R}$. We let $\widehat{q}: \mathcal{R}_v(^\rho\mathbb{R}) \to
\widehat{^\rho\mathbb{R}}$ be the quotient map.

For more details about $^\rho\mathbb{R}$ we refer the reader to A. Robinson \cite{aRob73}, A. H.
Lightstone and A. Robinson~\cite{LiRob}, W.A.J. Luxemburg~\cite{wLux}, and V. Pestov~\cite{vPes}. 

\section{The Quasi-Standard Part Mapping} \label{S: quasi}

We always view $^\rho\mathbb{R}$ as an extension of $\mathbb{R}$.

\begin{definition}\label{D:Feasible}{\em (Feasible):}
 Let $B \subseteq {^\rho\mathbb{R}}$. We will say that $B$ is \textbf{feasible} if 
$\mathbb{R}(B) \subseteq \mathcal{R}_v(^\rho\mathbb{R})$, where
$\mathbb{R}(B)$ denotes the subfield of $^\rho\mathbb{R}$ generated by $\mathbb{R}\cup B$. 
\end{definition}

The following proposition is immediate.

\begin{proposition} \label{P: Feasible} (a) For any $B \subseteq {^\rho\mathbb{R}}$, 
$B$ is feasible if and only if every non-zero element of
$\mathbb{R}(B)$ has valuation zero.

	(b) Suppose that $B \subseteq{^\rho\mathbb{R}}$ and, for any nonzero polynomial $P \in 
\mathbb{R}[t_1, t_2, ...,t_n]$ and any distinct $b_1, b_2, ..., b_n \in B$, $v(P(b_1, b_2, ..., b_n)) =
0$. Then $B$ is feasible. 

\end{proposition}

\begin{theorem} \label{T: Existence} {\em (Existence):} Let $B$ be feasible. Then:

	(a) There is a maximal field $\widehat{\mathbb{R}}$, which is a a subring of
$\mathcal{R}_v(^\rho\mathbb{R})$ that contains $\mathbb{R}(B)$, in symbols,
 $\mathbb{R}(B) \subset\widehat{\mathbb{R}}\subset \mathcal{R}_v(^\rho\mathbb{R})$

	(b) Any such field $\widehat{\mathbb{R}}$ is real-closed.

	(c) Any such field $\widehat{\mathbb{R}}$ is isomorphic to $\widehat{^\rho\mathbb{R}}$ 
via the quotient mapping $\widehat{q}$. 

\end{theorem}

\noindent{\bf Proof:} Results similar to this appear in various places in the
literature, such as Lemma 2 of MacLane \cite{sMacL}. We give a simple
independent proof.

(a) Given $B$, we obtain the desired field by applying Zorn's Lemma to
the collection of all fields $\mathbb{K}$ such that $\mathbb{R}(B) \subseteq \mathbb{K}
\subseteq \mathcal{R}_v(^\rho\mathbb{R})$, to obtain a maximal such field.

(b) Given a maximal field $\mathbb{K}$ as in (a), Let $\cl(\mathbb{K})$ be the relative
algebraic closure of $\mathbb{K}$ in $^\rho\mathbb{R}$. Since $^\rho\mathbb{R}$ is real-closed, so is
$\cl(\mathbb{K})$. So it suffices to show that $\cl(\mathbb{K}) \subset \mathcal{R}_v(^\rho\mathbb{R})$,
since this implies that $\cl(\mathbb{K}) = \mathbb{K}$ by the maximality of $\mathbb{K}$. So let 
$x \in \cl(\mathbb{K})$. Then $x$ is a root of some polynomial $t^n +
a_{n-1}t^{n-1} + ... + a_0$, with $a_k \in \mathbb{K}$. But the bound $|x| \le 1
+ |a_0| + ... + |a_{n-1}|$ implies that $v(x) \ge 0$, so $x \in\mathcal{R}_v(^\rho\mathbb{R})$.

(c) Since $\mathbb{K}\subseteq \mathcal{R}_v(^\rho\mathbb{R})$, the quotient mapping $\widehat{q}$,
restricted to $\mathbb{K}$, is a field homomorphism from $\mathbb{K}$ to $\widehat{^\rho\mathbb{R}}$.
Furthermore, $\widehat{q}$ is one-to-one, because if $\widehat{q}(x) =\widehat{q}(y)$,
then $v(x-y) > 0$, which implies that $x = y$ because $v$ is trivial on
$\mathbb{K}$. (I.e., $v(x) > 0$ for a nonzero $x \in \mathbb{K}$ would imply $v(1/x) < 0$,
contradicting that $\mathbb{K} \subseteq \mathcal{R}_v(^\rho\mathbb{R})$.) Thus $\widehat{q}$ defines
an isomorphism between $\mathbb{K}$ and its image $\widehat{q}(\mathbb{K})$. 
	It only remains to prove that $\widehat{q}(\mathbb{K}) = \widehat{^\rho\mathbb{R}}$. If not, 
then there's a 
$\xi \in \mathcal{R}_v(^\rho\mathbb{R})$ such that $\widehat{q}(\xi) \notin \widehat{q}(\mathbb{K})$.
Now let $P \in \mathbb{K}[t], P \ne 0$. We claim that $v(P(\xi)) = 0$.
Otherwise, we would have $v(P(\xi)) > 0$, since $v(\xi) = 0$. But then
$\widehat{q}(P(\xi)) = 0$, or equivalently, ${\widehat P}(\widehat{q}(\xi)) = 0$ in
$\widehat{^\rho\mathbb{R}}$, where $\widehat{P}$ is the polynomial $\widehat{q}(P)$ in
$\widehat{q}(\mathbb{K})[t]$. So $\widehat{q}(\xi)$ is algebraic over $\widehat{\mathbb{K}}$. But, 
since $\mathbb{K}$ is real-closed by (b), its isomorphic image $\widehat{q}(\mathbb{K})$ is a 
real-closed 
subfield of the ordered field $\widehat{^\rho\mathbb{R}}$. Thus $\widehat{q}(\xi) \in
\widehat{q}(\mathbb{K})$, contradicting the assumption about $\xi$. So $v(P(\xi)) = 0$. It follows that
$v(P(\xi)/Q(\xi)) = 0$ whenever 
$P,Q \in \mathbb{K}[t], P, Q \ne 0$. In other words, $\mathbb{K}(\xi) \subseteq
\mathcal{R}_v(^\rho\mathbb{R})$. But since
$\xi \notin \mathbb{K}$, this contradicts the maximality of $\mathbb{K}$. $\blacktriangle$

Note that this theorem is not vacuous, since we can let $B = \varnothing$.

Let $\widehat{j}: \widehat{^\rho\mathbb{R}} \to {^\rho\mathbb{R}}$ be the inverse of the
isomorphism $\widehat{q}$ described in Theorem~\ref{T: Existence}. So $\widehat{j}$ is an
isomorphic embedding of $\widehat{^\rho\mathbb{R}}$ in $^\rho\mathbb{R}$. We refer to $j$ as a 
\textbf{ $B$-invariant embedding} of $\widehat{^\rho\mathbb{R}}$ into
$^\rho\mathbb{R}$. Of course, this embedding is not unique, even for fixed
$B$.

	A field $\widehat{\mathbb{R}}$ that satisfies Theorem~\ref{T: Existence} is sometimes called a
\textbf{field of representatives} for $\mathcal{R}_v(^\rho\mathbb{R})$. We will also call it a
\textbf{$B$-copy of $\widehat{^\rho\mathbb{R}}$ in $^\rho\mathbb{R}$}. Henceforth, $\widehat{\mathbb{R}}$
always denotes such a field.

\begin{corollary} \label{C: DirectSum} 
$\mathcal{R}_v(^\rho\mathbb{R}) = \widehat{\mathbb{R}} \oplus \mathcal{I}_v(^\rho\mathbb{R})$. 
\end{corollary}

\Proof Let $\xi \in \mathcal{R}_v(^\rho\mathbb{R})$. By Theorem~\ref{T: Existence}, we know
there's an $x$ in $\widehat{\mathbb{R}}$ such that $\widehat{q}(\xi) = \widehat{q}(x)$. In other
words, $\xi = x + h$, for some $h \in \mathcal{I}_v(^\rho\mathbb{R})$. And since 
$v|\widehat{\mathbb{R}}$ is trivial, this representation $\xi = x + h$ is unique. $\blacktriangle$
\begin{definition} \label{D: Quasi} Given $\widehat{\mathbb{R}}$, we define the 
\textbf{quasi-standard part
mapping} $\widehat{\st}: \mathcal{R}_v(^\rho\mathbb{R}) \to \widehat{\mathbb{R}}$ by 
$\widehat{\st}(x+ h) = x$, for any $x \in \widehat{\mathbb{R}}$ and $h \in
\mathcal{I}_v(^\rho\mathbb{R})$. 
\end{definition}

	The following properties of $\widehat{\st}$ are straightforward:

\begin{proposition}[Properties] \label{P: Properties} 

	(a) $\widehat{\st}$ is an ordered ring homomorphism from $\mathcal{R}_v(^\rho\mathbb{R})$ onto
$\widehat{^\rho\mathbb{R}}$. 

	(b) $\widehat{\st}|\widehat{\mathbb{R}}$ is the identity. In particular, all real numbers and
all asymptotic numbers in $B$ are fixed points of $\widehat{\st}$.

	(c) $\widehat{\st} \circ \widehat{\st} = \widehat{\st}$.

	(d) The restriction of $\widehat{\st}$ to $\mathbb{R} \oplus \mathcal{I}_v(^\rho\mathbb{R})$ 
coincides
with $\st$, the usual standard part mapping in $^\rho\mathbb{R}$ (hence the
term ``quasi-standard part mapping''). 

	(e) $\widehat{\st} = {\widehat{j}} \circ \widehat{q}$. 
\end{proposition}

For applications in asymptotic analysis involving Laurent asymptotic
expansions or, more generally, expansions in Levi-Civita series (with
real or complex coefficients), it usually suffices to apply
Theorem~\ref{T: Existence} with $B = \varnothing$. For other applications, it can be
fruitful to use a non-empty $B$, in order to provide the mapping $\widehat{\st}$
with particular fixed points other than the reals. We conclude this
section by discussing a choice of $B$ that is useful for applications
involving asymptotic expansions with logarithm terms: 

\begin{notation}\label{N: Logs}{\em (Multiple Logarithms):} (a) Let $\lambda_1 = |{\ln} \rho|$
and, inductively, $\lambda_{n+1} = \ln(\lambda_n)$. Notice that the $\lambda_n$'s form a
decreasing sequence of infinitely large positive numbers in $^*\mathbb{R}$.

	(b) Let $l_n = q(\lambda_n)$, for $n \in \mathbb{N}$. So the $l_n$'s form a
decreasing sequence of infinitely large positive numbers in $^\rho\mathbb{R}$.
However, $v(l_n) = 0$. 
\end{notation}

Another way to define the $l_n$'s is to use the function
$^\rho{\ln}$  in $^\rho\mathbb{R}$, as discussed near the end of Section 1.
Specifically, $l_1 = {^\rho\ln}(1/s)$, and $l_{n+1} = {^\rho\ln}(l_n)$.

\begin{lemma} \label{L: Feasible} The set $\mathcal{L} = \{l_n : n \in \mathbb{R} \}$ is feasible. 
So are the sets $\mathcal{S}=\{e^{\pm i\pi/s}\ln^n{s}: n=0, 1, 2,\dots \}$ and $\mathcal{L}\cup
\mathcal{S}$. 
\end{lemma}

\Proof By Proposition~\ref{P: Feasible}, it suffices to prove that if $P
\in \mathbb{R}[t_1, t_2, ..., t_n]$ and $P \neq 0$, then $v(P(l_1, l_2, ...,
l_n)) = 0$. First note that this is obviously true when $P$ is a
monomial, since each $l_n$ has valuation $0$. If $P$ has more than one
term, we just need to show that $P(l_1, l_2, ..., l_n)$ has a ``dominant
term,'' compared to which all other terms are infinitesimal.
(Technically, this is the greatest term in the polynomial with respect
to the lexicographic ordering on the sequence of its exponents.) 
For this it suffices to show that  $m
> k$ implies $(\lambda_m)^p/(\lambda_k)^q
\approx 0$ whenever $p$ and $q$ are in $\mathbb{R}$. By the transfer principle, this is
equivalent to a limit statement in standard mathematics, which is easily
verified by l'Hopital's Rule. The sets $S$ and $\mathcal{L}\cup \mathcal{S}$ are treated similarly.
$\blacktriangle$ 

If we apply Theorem~\ref{T: Existence} with $B = \mathcal{L}\cup \mathcal{S}$, the resulting field of
representatives $\widehat{\mathbb{R}}$ contains all the $l_n$'s and $e^{\pm i\pi/s}\ln^n{s}$, which are
fixed points of $\widehat{\st}$. Under that assumption we have $\widehat{\st}(\ln{s}+r+s)=\ln{s}+r,\;
\widehat{\st}(e^{\pm i\pi/s}+r+s^2)=e^{\pm i\pi/s}+r$,\; $\widehat{\st}(s)=0$ where $r\in\mathbb{R}$.

\section{The Imbedding of\, $\widehat{^\rho\mathbb{R}}\bra t^{\mathbb{R}}\ket$ in $^\rho\mathbb{R}$}
\label{S: Imbedding}

Recall the definition in Section 1 of the field $\mathbb{R}\bra t^{\mathbb{R}}\ket$ We can
similarly define $\mathbb{K}\bra t^{\mathbb{R}}\ket$ for any field $\mathbb{K}$, and we will 
refer to all
fields of this type as \textbf{Levi-Civita fields}. A. Robinson
\cite{aRob73} showed that $\mathbb{R}\bra t^\mathbb{R})$ can be embedded in $^\rho\mathbb{R}$. We now
generalize Robinson's result (and Theorem~\ref{T: Existence} above) by providing an imbedding of
$\widehat{^\rho\mathbb{R}}\bra t^\mathbb{R}\ket$ in $^\rho\mathbb{R}$. Since $\widehat{^\rho\mathbb{R}}$ is
isomorphic to a subfield
$\widehat{\mathbb{R}}$  of $^\rho\mathbb{R}$, we temporarily identify $\widehat{^\rho\mathbb{R}}\bra
t^\mathbb{R}\ket$  with
$\widehat{\mathbb{R}}\bra t^\mathbb{R}\ket$

	For uniformity, we assume that any series $A \in {\widehat{^\rho\mathbb{R}}\bra t^\mathbb{R}\ket}$ is
written in the form
$A=\sum_{k=0}^\infty {{a_k} t^{r_k}}$, where the sequence $(r_k)$ of
reals is strictly increasing and unbounded, $a_k \in \widehat{^\rho\mathbb{R}}$ (or
$\widehat{\mathbb{R}}$), and if any $a_k = 0$, then $a_m = 0$ for every $m > k$.
Naturally, two such series are considered equal if any term that is in
one series but not the other has coefficient $0$. If the series  $\sum_{k=0}^\infty {a_k}\, x^{r_k}$ is
convergent in $^\rho\mathbb{R}$ for some $x \in {^\rho\mathbb{R}}$, we denote its sum by $A(x)$ and 
write
$A(x)=\sum_{k=0}^\infty {a_k}\, x^{r_k}$. Let $\mathcal{I}_v(^\rho\mathbb{R}_+) = \{h \in
{^\rho\mathbb{R}}_+ \mid v(h) > 0\}$ denote the multiplicative group of scales in $^\rho\mathbb{R}$, 
where
$^\rho\mathbb{R}_+$ is the set of positive elements of $^\rho\mathbb{R}$.
\begin{lemma} \label{L: Convergent} Let $A \in \widehat{^\rho\mathbb{R}}\bra t^\mathbb{R}\ket,\; h \in
\mathcal{I}_v(^\rho\mathbb{R}_+)$. Then

(a) The series $A(h)=\sum_{k=0}^\infty {a_k}\, h^{r_k}$ is convergent in $^\rho\mathbb{R}$.

(b) If $A(h) \ne 0$, then $v(A(h)) = r_0 v(h)$.
\end{lemma}

\Proof (a) If $A(h)$ is a finite series, there's nothing to prove.
Otherwise, since $v(a_k) \ge 0$ and $v(h) > 0$, it follows that $v({a_k}
h^{r_k}) \to \infty$ as $k \to \infty$. By remarks made
in Section 1, this implies that $a_k h^{r_k} \rightarrow 0$, which in
turn implies that $A(h)$ converges.

(b) This follows from the fact that $a_0 h^{r_0}$ is the dominant term
in any partial sum of $A(h)$, and the continuity of $v$. $\blacktriangle$

\begin{theorem}[The Imbedding] \label{T: Imbedding}
For any fixed $h \in \mathcal{I}_v(^\rho\mathbb{R}_+)$, the function 
$M_h: \widehat{^\rho\mathbb{R}}\bra t^\mathbb{R}\ket\to{^\rho\mathbb{R}}$ defined by $M_h(A) = A(h)$ is a
continuous ordered field embedding.
\end{theorem}

\Proof The previous lemma showed that $M_h$ is defined on 
$\widehat{^\rho\mathbb{R}}\bra t^\mathbb{R}\ket$.
The fact that $M_h$ is a ring homomorphism is clear. We also need
to show that $M_h$ is order-preserving: for any $A \in \widehat{^\rho\mathbb{R}}\bra t^\mathbb{R}\ket$, we
have $A > 0$ iff $a_0 > 0$ iff $A(h) > 0$, since $a_0 h^{r_0}$ is the
dominant term of $A(h)$. So $M_h$ is an ordered ring embedding. To show that $M_h$ is continuos,
note that $v(M_h(A)) = v(A(h)) = r_0 v(h)$ by the previous lemma, which equals $V(A)v(h)$, where
$V$ is the canonical valuation on $\widehat{^\rho\mathbb{R}}\bra t^\mathbb{R}\ket$. So $M_h$ is
valuation-preserving except for a factor of $v(h)$. This immediately implies continuity of $M_h$.
$\blacktriangle$

	Note that the embedding $M_h$ depends on the choice of $\widehat{\mathbb{R}}$ as well
as $h$ (not to mention $\rho$). The particular embedding $M_s$, where
$s=\overline{\rho}$ is the canonical scale of $^\rho\mathbb{R}$, is an extension of the embedding
$\Phi$, defined in A. Robinson
\cite{aRob73}. Unlike a typical $M_h$, $M_s$ is also valuation-preserving.

\begin{corollary} \label{C: Algorithm} Let $\widehat{^\rho\mathbb{R}}\bra h^\mathbb{R}\ket$ denote the
range of $M_h$. The inverse of the mapping $M_h$ is given explicitly as follows: for any $y$ in
$\widehat{^\rho\mathbb{R}}\bra h^\mathbb{R}\ket$, $y = M_h(A) = A(h)$, where the series
$A=\sum_{k=0}^\infty {a_k}\, t^{r_k}$ is defined inductively by:
\begin{align}
r_n = v \left(y - \sum_{k<n}{a_k h^{r_k}} \right)/ v(h),\quad
a_n = \widehat{\st}\left( (y -\sum_{k<n}{a_k h^{r_k}})/h^{r_n} \right),\notag
\end{align}
with the understanding that if this algorithm ever gives $r_n
= \infty$, the series terminates at the previous term. 
\end{corollary}

\Proof We verify the algorithm for $r_0$ and $a_0$. Its validity
for subsequent terms follows by the same reasoning. First of all, if\, $y = 0$,
then the algorithm gives $r_0 = \infty$, whence $A$ is the ``null series,'' which is the zero
element of $\widehat{^\rho\mathbb{R}}\bra t^\mathbb{R}\ket$. If $y \neq 0$, we know from Lemma~\ref{L:
Convergent}(b) that
$v(y) = r_0v(h)$, so $r_0 = v(y)/v(h)$, as desired. It then follows that $v(y/h^{r_0}) = v(y) -
r_0v(h) = 0$, so $y/h^{r_0} \in \mathcal{R}_v(^\rho\mathbb{R})$. Therefore, $\widehat{\st}(y/h^{r_0})$ is
well-defined. Furthermore, by definition of $\widehat{\st}$, the value of $a_0$ defined by our algorithm is
the unique member of $\widehat{^\rho\mathbb{R}}$ (more precisely, of the
field of representatives $\widehat{\mathbb{R}}$) such that
$v(y/h^{r_0} -a_0) > 0$. But the correct $a_0$ must satisfy $v(y - a_0h^{r_0}) = v(a_1h^{r_1}) =
r_1v(h) > r_0v(h))$, or, equivalently, $v(y/h^{r_0} - a_0) > 0$, so the defined value of $a_0$ is
correct. $\blacktriangle$ 

\section{The Hahn Field Representation of $^\rho\mathbb{R}$} \label{S: Hahn}

In this section we prove the existence of valuation-preserving
isomorphisms between $^\rho\mathbb{R}$ and the Hahn field $\widehat{^\rho\mathbb{R}}(t^\mathbb{R})$. A
similar result in the setting of ultraproducts appears in B. Diarra \cite{bDi} (Corollary
to Proposition 8). Our result may be more general in a couple of
ways. In particular, our isomorphisms allows us to specify fixed points,
as discussed at the end of Section 2.

Unless stated otherwise, the results in this section pertain to general
Krull's valuations, not necessarily convex or real-valued (W. Krull~\cite{wKrull}).

\begin{definition} \label{D: Extensions} {\em (Maximal Extensions):}
Let $(\mathbb{K}_i, v_i)$ be valuation fields, $i = 1, 2$. If
$\mathbb{K}_1$ is a subfield of $\mathbb{K}_2$ and $v_2 | \mathbb{K}_1 = v_1$, we say that 
$(\mathbb{K}_2, v_2)$ is an {\bf extension} of $(\mathbb{K}_1, v_1)$. If, in addition, the value
groups $G_{v_1}$ and $G_{v_2}$ are equal and the residue class fields
$\widehat{\mathbb{K}}_1$ and $\widehat{\mathbb{K}}_2$ coincide (in the sense that the
natural embedding of $\widehat{\mathbb{K}}_1$ in $\widehat{\mathbb{K}}_2$ is
surjective), the extension is called \textbf{immediate}. A valuation field
with no proper immediate extensions is called \textbf{maximal}. 
\end{definition}

W.A.J. Luxemburg \cite{wLux} proved that maximality is equivalent to spherical
completeness in fields with convex, real-valued valuations.
W. Krull \cite{wKrull} proved that every valuation field has a maximal
immediate extension. We are more interested in the following result, due
to Kaplansky \cite{iKa} (characteristic 0 case of Theorem 5):

\begin{theorem} \label{T: Kaplansky} {\em (Kaplansky):} Let $(\mathbb{K}, v)$ be a valuation field whose
residue class field has characteristic 0. Then the maximal immediate extension
of $(\mathbb{K}, v)$ is unique, in the sense that between any two maximal
immediate extensions of $(\mathbb{K}, v)$ there is a valuation-preserving
isomorphism that is the identity on $\mathbb{K}$. 
\end{theorem}

\begin{corollary} \label{C: Extension} {\em (Extension):}
Let $(\mathbb{K}_i, v_i)$, for $i = 1,2$, be a valuation
field whose residue class field has characteristic 0, and let $(\mathbb{L}_i,
w_i)$ be a maximal immediate extension of $(\mathbb{K}_i, v_i)$. Then any
valuation-preserving isomorphism between $\mathbb{K}_1$ and $\mathbb{K}_2$ can be extended
to such an isomorphism between $\mathbb{L}_1$ and\, $\mathbb{L}_2$. 

\end{corollary}

\begin{theorem} \label{T: Isomorphism} {\em (Isomorphism):} 
For any choice of the field of representatives $\widehat{\mathbb{R}}$ of $\widehat{^\rho\mathbb{R}}$, 
there
is a valuation-preserving isomorphism
$J: {^\rho\mathbb{R}}\to\widehat{^\rho\mathbb{R}}(t^\mathbb{R})$ with the following additional
properties:

(a) $J^{-1}$ restricted to $\widehat{^\rho\mathbb{R}}\bra t^\mathbb{R}\ket$ is the mapping 
$M_s$ defined in Section~\ref{S: Imbedding}.

(b) $J$ maps the scale $s$ of $^\rho\mathbb{R}$ to the indeterminate\, $t$\, of
$\widehat{^\rho\mathbb{R}}(t^\mathbb{R})$.

(c) $J$ restricted to $\widehat{^\rho\mathbb{R}}$ is the quotient mapping\, $\widehat{q}$.

\end{theorem}

\Proof We apply Corollary~\ref{C: Extension} with $\mathbb{K}_1 = \widehat{^\rho\mathbb{R}}\langle
s^\mathbb{R}\rangle$, $\mathbb{K}_2 = \widehat{^\rho\mathbb{R}}(t^\mathbb{R})$, 
$\mathbb{L}_1 = {^\rho\mathbb{R}}$, and
$\mathbb{L}_2 = 
\widehat{^\rho\mathbb{R}}(t^\mathbb{R})$ (all with their usual valuation), and the valuation-preserving
isomorphism
$({M_s})^{-1}$ between $\mathbb{K}_1$ and $\mathbb{K}_2$ (see Section~\ref{S: Imbedding}).
$\mathbb{L}_1$ is an immediate extension of $\mathbb{K}_1$ because both have valuation
group $\mathbb{R}$ and residue class field $\widehat{^\rho\mathbb{R}}$ (and the residue class
fields coincide in the sense described in the definition of
immediateness). Similarly,\, $\mathbb{L}_1$ is an immediate extension of $\mathbb{K}_1$.
The maximality of $(^\rho\mathbb{R}, v)$ is proved in W.A.J. Luxemburg~\cite{wLux}. The
maximality of Hahn fields is proved in W. Krull~\cite{wKrull}. Thus, by
Corollary~\ref{C: Extension}, there is a valuation-preserving isomorphism $J$
between $^\rho\mathbb{R}$ and $\widehat{^\rho\mathbb{R}}(t^\mathbb{R})$ that satisfies (a). 
Properties (b)
and (c) follow directly from the definition of $M_s$. $\blacktriangle$ 

\begin{remark} Note that the isomorphism provided by Theorem~\ref{T: Isomorphism} is
non\-unique and nonconstructive at several stages. First, we have some
freedom in the choice of the set $B$. Then, the definition of $^\rho\mathbb{R}$
(Theorem~\ref{T: Existence}) uses Zorn's Lemma. Finally, the definition of the
isomorphism in Kaplansky's result (Theorem~\ref{T: Kaplansky}) is also
nonconstructive. We could also use any $M_h$ instead of $M_s$ as the
basis for $J$ if we don't need $J$ to be valuation-preserving.
\end{remark}
\begin{remark} Theorem~\ref{T: Isomorphism} and most of our other results so far about $^\rho\mathbb{R}$
hold for any spherically complete, real-closed extension of $\mathbb{R}$ with a convex, real-valued
nonarchimedean valuation.
\end{remark}
\section{Isomorphism Between $^*\mathbb{R}$ and $\widehat{^\rho\mathbb{R}}$}\label{S: Isomormorphism}

Given $^*\mathbb{R}$ one might ask if there is essentially just one field
$^\rho\mathbb{R}$. In other words, does the structure of $^\rho\mathbb{R}$ depend on the
choice of $\rho$? This question leads to two cases. For positive
infinitesimals $\rho_1$ and $\rho_2$, suppose that $\ln\rho_1/\ln\rho_2$ is finite but
not infinitesimal. This condition is easily seen to be equivalent to ${\cal{M}}_{\rho_1}(^*\mathbb{R})
=  {\cal{M}}_{\rho_2}(^*\mathbb{R})$. Therefore, since the maximal ideal 
$\mathcal{N}_\rho(^*\mathbb{R})$
is determined uniquely by $\mathcal{M}_\rho(^*\mathbb{R})$, this condition implies that
$^{\rho_1}\!\mathbb{R}$ and $^{\rho_2}\!\mathbb{R}$ are not just isomorphic
but identical. Without assuming this condition, we can still prove that any two
$^\rho\mathbb{R}$'s defined from a particular $^*\mathbb{R}$ are isomorphic, but only
under the assumption that $^*\mathbb{R}$ is a special model (defined later in
this section). In fact, under this assumption, $\widehat{^\rho\mathbb{R}}$ (for any
$\rho$) is actually isomorphic to $^*\mathbb{R}$. This is perhaps surprising,
since it seems that $^\rho\mathbb{R}$ is a ``much smaller'' field than $^*\mathbb{R}$,
and $\widehat{^\rho\mathbb{R}}$ in turn is ``much smaller'' than $^\rho\mathbb{R}$.

$^*\mathbb{R}$ is called \textbf{fully saturated} if it is
$\card(^*\mathbb{R}$)-saturated. We will first prove our results for fully
saturated $^*\mathbb{R}$, and then generalize them to special $^*\mathbb{R}$.

We briefly review some basic concepts of model theory: Let $\mathcal{L}$ be
a first-order language. (For our purposes, $\mathcal{L}$ is the first-order
language of an ordered ring, with two binary operation symbols and one
binary relation symbol besides equality). Two $\mathcal{L}$-structures
$\mathcal{A}$ and $\mathcal{B}$ are called {\bf elementarily equivalent},
written $\mathcal{A} \equiv \mathcal{B}$, if $\mathcal{A}$ and $\mathcal{B}$ satisfy
exactly the same statements of $\mathcal{L}$ (with no free variables).
Also, $\mathcal{A}$ is said to be \textbf{elementarily embedded} in $\mathcal{B}$, 
written $\mathcal{A} \prec{\mathcal{B}}$, if $\mathcal{A}$ is a substructure of $\mathcal{B}$ and
$\mathcal{A}$ and
$\mathcal{B}$ satisfy exactly the same formulas of $\mathcal{L}$ with free variables interpreted as any
members of the domain of $\mathcal{A}$.

The second condition in the definition of $\mathcal{A} \prec {\mathcal{B}}$ is
much stronger than $\mathcal{A} \equiv \mathcal{B}$. For further details on
these concepts, see Chang and Keisler~\cite{CKeis}, which also includes the
following important results:

\begin{theorem}[Model Theory of Fields] \label{T: Model Theory} 

	(a) {\em Tarski Theorem}: All real-closed fields are elementarily
equivalent (Theorem 5.4.4 in \cite{CKeis}).

	(b) {\em Model-Completeness:} If $\mathcal{B}$ is a real-closed field, then $\mathcal{A} \prec
\mathcal{B}$ if and only if $\mathcal{A}$ is a real-closed subfield of $\mathcal{B}$ (page 110 in
\cite{CKeis}).

	(c) {\em Elementary Chain Theorem:} Let $\{\,{\mathcal{A}}_i \mid i \in I \,\}$ be a
family of structures such that $I$ is well-ordered, and ${\mathcal{A}}_i
\prec {\mathcal{A}}_j$ whenever $i, j \in I$ and $i < j$. Then ${\mathcal{A}}_i
\prec \cup_{i \in I} {\cal{A}}_i$, for every $i \in I$ (Theorem 3.1.13 in \cite{CKeis}).

(d) {\em Uniqueness of special and fully saturated models:} Any two
elementarily equivalent, special models of the same cardinality are
isomorphic (and every fully saturated model is special). (Theorems
5.1.17 and 5.1.13. in \cite{CKeis}) 
\end{theorem}
\begin{theorem} \label{T: Fully Saturated} If {\em $^*\mathbb{R}$} is fully saturated, then
$^*\mathbb{R}$ is isomorphic to any residue class field $\widehat{^\rho\mathbb{R}}$ as defined in 
Section~\ref{S: Asy}. 
\end{theorem}

\Proof Assume $^*\mathbb{R}$ is fully saturated, and let $\kappa = \card(^*\mathbb{R})$. We will use
Theorem~\ref{T: Model Theory}~(d). $^*\mathbb{R}$ is real-closed by the transfer principle, and
$\widehat{^\rho\mathbb{R}}$ is real-closed by Theorem~\ref{T: Existence}, (b) and (c). So, by
Theorem~\ref{T: Model Theory} (a), $^*\mathbb{R}
\equiv {^\rho\mathbb{R}}$. Recall that $^\rho\mathbb{R}$ is a quotient field of a subring of 
$^*\mathbb{R}$,
and
$\widehat{^\rho\mathbb{R}}$ is defined similarly from $^\rho\mathbb{R}$. We now combine these
operations: let $Q$ be the composition $\widehat{q} \circ q$. The domain of
$Q$ is 
\[
\mathcal{F}_\rho(^*\mathbb{R})=\{x\in {^*\mathbb{R}}\mid |x|<1/\sqrt[n]{\rho}\; \text{for all}\;
n\in\mathbb{N}\}.
\] 

The kernel of $Q$ is 
\[
\mathcal{I}_\rho(^*\mathbb{R})=\{x\in^*\mathbb{R} \mid |x|<\sqrt[n]{\rho}\; \text{for
some}\; n\in\mathbb{N}\}.
\] 
It is clear from the above that $\card(\widehat{^\rho\mathbb{R}}) \le \kappa$. Therefore, 
if we can show
that
$\widehat{^\rho\mathbb{R}}$ is $\kappa$-saturated, we are done because $\widehat{^\rho\mathbb{R}}$ 
is then
fully saturated and has the same cardinality as $^*\mathbb{R}$. 
(Here we are using the simple fact that a
$\kappa$- saturated model must have cardinality at least $\kappa$.) So let 
\[
\mathcal{C} = \{ (a_i, b_i) | i\in I\},
\]
be a collection of fewer than $\kappa$ open intervals in $\widehat{^\rho\mathbb{R}}$ with the F.I.P. 
Let $A =\{a_i\}$ and $B = \{b_i\}$. For each $i \in I$, let $\tilde a_i \in
{^*\mathbb{R}}$ be such that $Q(\tilde a_i) = a_i$. Similarly let $Q(\tilde b_i)
= b_i$ for every $i \in I$. Our goal is to prove that $\cal{C}$ has nonempty intersection. We
proceed by cases. If $A$ has a greatest element and $B$ has a least
element, this follows trivially. If $A$ has no greatest element and $B$ has no least element, let
$\tilde{\cal C} = \{ (\tilde a_i, \tilde b_i) | i \in I\}$. Since $\cal{C}$
has the FIP, so does $\tilde{\cal{C}}$. Since $^*\mathbb{R}$ is fully saturated, there
is some $x \in \bigcap \tilde{\mathcal{C}}$. It is then easily seen that $Q(x)
\in \bigcap \mathcal{C}$.Finally, we consider the case that $A$ has a greatest element $a$, but
$B$ has no least element. (The reverse situation is handled similarly).
Now let $\tilde a$ be such that $Q(\tilde a) = a$, and define\, $\tilde{\mathcal C} = \{ (\tilde a +
\rho^{1/n}, \tilde b_i) | i \in I, n \in \mathbb{N}\}$. $\tilde{\mathcal{C}}$ has cardinality less than
$\kappa$ and has the F.I.P., so $\bigcap \tilde{\mathcal{C}} \neq \varnothing$. Let $x \in
\bigcap \tilde{\mathcal{C}}$. Clearly, $Q(x) < b_i$ for every $i \in I$. And
since $x > \tilde a + \rho^{1/n}$ for every $n$, $Q(x) > Q(\tilde a) =
a$. Thus $Q(x) \in \bigcap\mathcal{C}$. $\blacktriangle$

The following corollary is immediate.

\begin{corollary} \label{C: Two Rho's} Let $\rho_1$ and $\rho_2$ be positive
infinitesimals in the same fully saturated {\em $^*\mathbb{R}$}. Then:

	(a) The residue class fields\, $\widehat{^{\rho_1}\!{\mathbb{R}}}$ and\,
 $\widehat{^{\rho_2}\!{\mathbb{R}}}$ are isomorphic.

	(b) The fields $^{\rho_1}\!{\mathbb{R}}$ and\, $^{\rho_2}\!{\mathbb{R}}$  are isomorphic.

	(c) In fact, any $^\rho\mathbb{R}$ is isomorphic to the Hahn field {\em $^*\mathbb{R}(t^{\mathbb{R}})$}. 
\end{corollary}
\begin{remark} 
Parts (a) and (b) of this corollary appear on page 196 of W.A.J. Luxemburg~\cite{wLux}, without
proof and without stating the necessity of any saturation assumption beyond
$\aleph_1$-saturation. 
\end{remark}

	Since we are viewing $^*\mathbb{R}$, $^\rho\mathbb{R}$, and $\widehat{^\rho\mathbb{R}}$ as extension
fields of $\mathbb{R}$, it is desirable to know when these isomorphisms can
be specified to be ``over $\mathbb{R}$'':
\begin{corollary} \label{C: Fixed R} Assume that $^*\mathbb{R}$ is fully saturated, and
$\card(^*\mathbb{R}) >  \card(\mathbb{R})$. Then the isomorphisms described in
Theorem~\ref{T: Fully Saturated} and Corollary~\ref{C: Two Rho's} can be chosen to be the
identity on $\mathbb{R}$. 
\end{corollary}
\Proof We know that $\mathbb{R}$ is a real-closed subfield of both the
real-closed fields $^*\mathbb{R}$ and $^\rho\mathbb{R}$. So, by Theorem~\ref{T: Model Theory}~(b),
$\mathbb{R} \prec {^*\mathbb{R}}$ and $\mathbb{R} \prec \widehat{^\rho\mathbb{R}}$. These conditions are
clearly equivalent to $\mathbb{R}^\prime \equiv {^*\mathbb{R}}^\prime$ and $\mathbb{R}^\prime \equiv
\widehat{^\rho\mathbb{R}}^\prime$, where these three new structures are obtained by adding
individual constants for every real number to $\mathbb{R}$, $^*\mathbb{R}$, and
$^\rho\mathbb{R}$, respectively. (These new structures are $\mathcal{L}^\prime$-structures, where
$\mathcal{L}^\prime$ is the first-order language of an ordered ring, augmented by a constant 
symbol for each
real number.) Thus, $^*\mathbb{R}^\prime \equiv {^\rho\mathbb{R}^\prime}$. Also, since 
$^*\mathbb{R}$ and $\widehat{^\rho\mathbb{R}}$ are $\card(\mathbb{R})^+$-saturated, augmenting them by 
adding $\card(\mathbb{R})$ constants does not decrease their saturation (C.C. Chang and H.J.
Keisler~\cite{CKeis}, Proposition 5.1.1(iii)). In other words, $^*\mathbb{R}^\prime$ and
$^\rho\mathbb{R}^\prime$ are also fully saturated. Therefore, by Theorem~\ref{T: Model Theory}~(d),
$^*\mathbb{R}^\prime$ and $\widehat{^\rho\mathbb{R}}^\prime$ are isomorphic. But an isomorphism between
$^*\mathbb{R}^\prime$ and
$^\rho\mathbb{R}^\prime$ is simply an isomorphism between $^*\mathbb{R}$ and $\widehat{^\rho\mathbb{R}}$
that is the identity on $\mathbb{R}$. $\blacktriangle$

	It is well known that fully saturated models are hard to come by. To
prove the existence of a fully saturated $^*\mathbb{R}$ requires assuming
either some form of the continuum hypothesis, or that $\card(^*\mathbb{R})$ is
inaccessible. For example, if $^*\mathbb{R}$ is a countable ultrapower of
$\mathbb{R}$, then $^*\mathbb{R}$ is $\aleph_1$-saturated and $\card(^*\mathbb{R}) =
2^{\aleph_0}$. So the ordinary continuum hypothesis would make this
$^*\mathbb{R}$ fully saturated. But without a special assumption such as the continuum hypothesis, we
cannot prove that Theorem~\ref{T: Fully Saturated} and Corollaries~\ref{C: Two Rho's} and
\ref{C: Fixed R} are non-vacuous. We will now fix this problem by considering
special models.
\begin{definition} Let $\mathcal{A}$ be a structure with domain $A$. $\cal A$
is called \textbf{special} if $\mathcal{A}$ is the union of a chain of
structures $\{ {\cal A}_\kappa \}$, where $\kappa$ ranges over all
cardinals less than {\em Card}($A$), each ${\cal A}_\kappa$ is
$\kappa^+$-saturated, and ${\cal A}_{\kappa_1} \prec {\cal
A}_{\kappa_2}$ whenever $\kappa_1 < \kappa_2$. The family $\{ {\cal
A}_\kappa \}$ is called a {\bf specializing chain} for $\cal A$.
\end{definition}

		For structures whose cardinality is regular, specialness is equivalent
to full saturation. However, while fully saturated models of
singular cardinality cannot exist, the existence of a special $^*\mathbb{R}$
(of an appropriate singular cardinality) can be proved in ZFC (see
Jin~\cite{rJ}, Proposition 3.4).
\begin{theorem} \label{T: Special} If $^*\mathbb{R}$ is special (instead of fully
saturated), the conclusions of Theorem~\ref{T: Fully Saturated} and Corollaries~\ref{C: Two
Rho's} and \ref{C: Fixed R} still hold. 
\end{theorem}
\Proof Assume $^*\mathbb{R}$ is special, and let $\mu = {\rm
\card}(^*\mathbb{R})$. Let $\{ {\cal A}_\kappa : \kappa < \mu \}$ be a
specializing chain for $^*\mathbb{R}$, and let any positive infinitesimal
$\rho$ be given. Since $^*\mathbb{R}$ is the union of the $\mathcal{A}_\kappa$'s,
we can choose $\nu < \mu$ such that $\rho \in {\cal A}_\nu$. Now replace
every ${\cal A}_\kappa$ with $\kappa < \nu$ by this ${\cal A}_\nu$. In
this way we get a new specializing chain for $^*\mathbb{R}$\, such that $\rho
\in {\cal A}_\kappa$ for every $\kappa < \mu$. So, for every $\kappa < \mu$, we can
construct quotient structures from ${\cal A}_\kappa$ that are analogous to $^\rho\mathbb{R}$ and 
$\widehat{^\rho\mathbb{R}}$,
with the constructions combined into one step as in the proof of
Theorem~\ref{T: Fully Saturated}. That is, let 
\[
D_\kappa = \{x \in {\cal A}_\kappa : |x|
< \rho^{-1/n}, {\rm for \; every} \; n \in \mathbb{N}\},
\]
and let 
\[
K_\kappa = \{x \in {\cal A}_\kappa : |x| < \rho^{1/n}, {\rm for \; some} \; n \in \mathbb{N}\}.
\]
As before, $K_\kappa$ is the maximal convex ideal of the
convex ring $D_\kappa$. Let ${\cal B}_\kappa$ be the quotient field
$D_\kappa/K_\kappa$, and let $Q_\kappa$ be the quotient map. There is a canonical embedding of
each ${\cal B}_\kappa$ in $\widehat{^\rho\mathbb{R}}$. Identifying each ${\cal B}_\kappa$ with its image in
$\widehat{^\rho\mathbb{R}}$, it is simple to show that $\widehat{^\rho\mathbb{R}}$ is the union of the
increasing sequence of subfields $\{ {\cal B}_\kappa : \kappa < \mu \}$. Our next goal is to
show that these subfields of $\widehat{^\rho\mathbb{R}}$ form a specializing chain for
$\widehat{^\rho\mathbb{R}}$. By Theorem~\ref{T: Model Theory} (c), each ${\cal A}_\kappa \prec {^*\mathbb{R}}$. Therefore,
since $^*\mathbb{R}$ is real-closed, so is each ${\cal A}_\kappa$. It follows,
as in the proof of Theorem~\ref{T: Existence}, (b) and (c), that each ${\cal
B}_\kappa$ is real-closed. Thus, if $\kappa_1 < \kappa_2 < \mu$, the
fact that ${\cal B}_{\kappa_1} \subseteq {\cal B}_{\kappa_2}$ implies
that ${\cal B}_{\kappa_1} \prec {\cal B}_{\kappa_2}$. Also, since each ${\cal A}_\kappa$ is
$\kappa^+$-saturated, the proof of Theorem~\ref{T: Fully Saturated} shows that each ${\cal
B}_\kappa$ is $\kappa^+$-saturated as well. (I.e., the argument there that shows that
any level of saturation satisfied by $^*\mathbb{R}$ is passed on to $\widehat{^\rho\mathbb{R}}$; full
saturation is not required for this.) So we have shown that $\{ {\cal B}_\kappa : \kappa < \mu \}$ is a
specializing chain for $\widehat{^\rho\mathbb{R}}$, and therefore $\widehat{^\rho\mathbb{R}}$ is special.
We know that $\card(\widehat{^\rho\mathbb{R}}) \leq \mu = \card(^*\mathbb{R})$. But for each cardinal
$\kappa < \mu$,
$\widehat{^\rho\mathbb{R}}$ contains the subfield ${\cal B}_\kappa$, which is $\kappa^+$-saturated and
therefore has cardinality greater than $\kappa$. Therefore, $\card(\widehat{^\rho\mathbb{R}}) =
\card(^*\mathbb{R})$. As in the proof of Theorem~\ref{T: Fully Saturated}, $^*\mathbb{R}$ and 
$\widehat{^\rho\mathbb{R}}$ are
elementarily equivalent. Therefore, by Theorem~\ref{T: Model Theory}~(d), $^*\mathbb{R}$ and 
$\widehat{^\rho\mathbb{R}}$ are
isomorphic. The conclusions of Corollaries~\ref{C: Two Rho's} and \ref{C: Fixed R} follow immediately.
$\blacktriangle$

		Here is a brief summary of the main results of this section:
\begin{corollary} \label{to} (a) There exists a {\em $^*\mathbb{R}$} that satisfies the
conclusions of Theorem~\ref{T: Fully Saturated} and Corollaries~\ref{C: Two Rho's} and \ref{C:
Fixed R}.

(b) Assume GCH. Then, for every regular uncountable cardinal $\kappa$,
there is a {\em $^*\mathbb{R}$} of cardinality $\kappa$ that satisfies the
conclusions of Theorem~\ref{T: Fully Saturated} and Corollary~\ref{C: Two Rho's}, and, if $\kappa
>$ \card$(\mathbb{R})$, Corollary \ref{C: Fixed R}. 
\end{corollary}

\Proof (a) Immediate by Theorem~\ref{T: Special} and the remark preceding it.

(b) For any infinite cardinal $\kappa$, it is possible (using a
so-called $\kappa^+$-good ultrafilter) to construct a $^*\mathbb{R}$ that is
$\kappa^+$-saturated and has cardinality $2^\kappa$. (See T. Lindstrom
\cite{tLin}, Theorem III.1.3.) So GCH implies that there exists a fully
saturated $^*\mathbb{R}$ of each uncountable regular cardinality. $\blacktriangle$ 

{\bf Acknowledgment:} The authors wish to thank Lou van den Dries, Ward
Henson, Jerome Keisler, and David Marker for helpful discussions of the
preliminary version of this paper.

\end{document}